\documentclass[a4paper,11pt,oneside,reqno]{amsart}
\usepackage{amsmath}
\usepackage{amsthm}
\usepackage{amssymb}
\usepackage{graphicx}
\usepackage{mathrsfs}
\theoremstyle{plain}
\newtheorem{teore}{Theorem}[section]
\newtheorem{defin}[teore]{Definition}
\newtheorem{lem}[teore]{Lemma}
\newtheorem{coro}[teore]{Corollary}
\newtheorem{propo}[teore]{Proposition}
\newtheorem*{claim}{Claim}
\newtheorem*{claim*}{Claim}
\theoremstyle{remark}
\newtheorem{ejemplo}[teore]{{\sc Example}}
\newtheorem{notas}[teore]{{\sc Remark}}
\newcommand{\nrm}[1]{\|#1\|}

\newcommand{\prop}{\begin{propo}}
\newcommand{\fprop}{\end{propo}}
\newcommand{\cor}{\begin{coro}}
\newcommand{\fcor}{\end{coro}}

\newcommand{\defi}{\begin{defin}\rm}
\newcommand{\fdefi}{\end{defin}}
\newcommand{\eje}{\begin{ejemplo}}
\newcommand{\feje}{\end{ejemplo}}
\newcommand{\lema}{\begin{lem}}
\newcommand{\flema}{\end{lem}}
\newcommand{\teor}{\begin{teore}}
\newcommand{\fteor}{\end{teore}}
\newcommand{\nota}{\begin{notas}\rm}
\newcommand{\fnota}{ \end{notas}}
\newcommand{\clam}{\begin{claim}}
\newcommand{\fclam}{\end{claim}}
\newcommand{\clams}{\begin{claim*}}
\newcommand{\fclams}{\end{claim*}}

\newcommand{\lclam}{\begin{lclaim}}
\newcommand{\flclam}{\end{lclaim}}
\newcommand{\prucl}{\prue[Proof of Claim:]}
\newcommand{\fprucl}{\fprue}
\newcommand{\ben}{\begin{enumerate}}
\newcommand{\een}{\end{enumerate}}
\newcommand{\bit}{\begin{itemize}}
\newcommand{\eit}{\end{itemize}}

\newcommand{\mc}[1]{\mathcal{#1}}

\newcommand{\mk}[1]{\mathfrak{#1}}

\newcommand{\casos}{\begin{itemize}}
\newcommand{\fcasos}{\end{itemize}\setcounter{cs}{1}}

\newcommand{\fin}{\textsc{FIN}}

\newcommand{\conj}[2]{ \{ {#1}\,:\,{#2} \} }

\newcommand{\om}{\omega}

\newcommand{\ip}{\sqsubseteq}

\newcommand{\buit}{\emptyset}
\newcommand{\0}{\emptyset}

\newcommand{\ga}{\gamma}

\newcommand{\al}{\alpha}
\newcommand{\be}{\beta}
\newcommand{\de}{\delta}
\newcommand{\De}{\Delta}

\newcommand{\sig}{\sigma}

\newcommand{\vphi}{\varphi}

\newcommand{\vep}{\varepsilon}

\newcommand{\eqs}{{\mathfrak X}}

\newcommand{\N}{{\mathbb N}}

\newcommand{\rest}{\upharpoonright}

\newcommand{\supp}{\mathrm{supp\, }}

\newcommand{\con}{\subseteq}

\newcommand{\prue}{\begin{proof}}
\newcommand{\fprue}{\end{proof}}
\makeindex

\begin{document}
\title[Pre-compact families of finite sets and weakly null sequences]{Pre-compact families of finite
sets of integers and weakly null sequences in Banach spaces}

\author{J. Lopez-Abad}
\address{Equipe de Logique Math\'{e}matique, Universit\'{e} Paris 7, 2 place Jussieu, 75251 Paris Cedex 05, France.}
\email{abad@logique.jussieu.fr}
\author{S. Todorcevic}
\address{C.N.R.S., U.M.R. 7056 et Universit\'{e} Paris 7, U.F.R. de Math\'{e}matiques, Case 7012, 2 Place Jussieu,
75251 Paris, Cedex 5, France. \newline Department of Mathematics, University of Toronto, Toronto, Canada, M5S
3G3}
 \email{stevo@math.toronto.edu}

\subjclass[2000]{Primary 05D10, 46B20, 46A35;
  Secondary 03E05}%
\maketitle

\section{Introduction}
The affinities between the infinite-dimensional Ramsey theory and
some problems of the Banach space theory and especially those
dealing with Schauder basic sequences have been explored for quite
some time, starting perhaps with Farahat's proof of Rosenthal's
$\ell_1$-theorem (see \cite{li-tza} and \cite{odell}). The
Nash-Williams' theory though implicit in all this was not fully
exploited in this context. In this paper we  try to demonstrate
the usefulness of this theory by applying it to the classical
problem of finding unconditional basic-subsequence of a given
normalized weakly null sequence in some Banach space $E$. Recall
that Bessaga and Pelczynski \cite{bess-pel} have shown that every
normalized weakly null sequence in a Banach space contains a
subsequence forming a Schauder basis for its closed linear span.
However, as demonstrated by Maurey and Rosenthal \cite{mau-ros}
there exist weakly null sequences in Banach spaces without
unconditional basic subsequences. So one is left with a task of
finding additional conditions on a given weakly null sequence
guaranteeing the existence of unconditional subsequences. One such
condition, given by Rosenthal himself around the time of
publication of \cite{mau-ros} (see also \cite{odell}). When put in
a proper context Rosenthal's condition reveals the connection with
the Nash-Williams theory. It says that if a weakly null sequence
$(x_n)$ in some space of the form $\ell_{\infty}(\Gamma)$ is such
that each $x_n$ takes only the values $0$ or $1$, then $(x_n)$ has
an unconditional subsequence. To see the connection, consider the
family
$$\mc F = \{\{n\in \N: x_n(\gamma)=1\}: \gamma\in \Gamma \}$$
and note that $\mc F$ is a pre-compact family of finite subsets of
$\N.$ As pointed out in \cite{odell}, Rosenthal result is
equivalent saying that there is an infinite subset $M$ of $\N$
such that the {\it trace} $${\mc F}[M]=\{t\cap M: t\in {\mc F}
\}$$ is {\it hereditary}, i.e., it is downwards closed under
inclusion. On the other hand, recall that the basic notion of the
Nash-Williams' theory is the notion of a {\it barrier}, which is
simply a family $\mc F$ of finite subsets of $\N$ no two members
of which are related under the inclusion which has the property
that an arbitrary infinite subset of $\N$ contains an initial
segment in $\mc F.$ Thus, in particular, $\mc F$ is a pre-compact
family of finite subsets of $\N.$ Though the trace of an arbitrary
pre-compact family might be hard to visualize, a trace ${\mc
B}[M]$ of a barrier $\mc B$ is easily to compute as it is simply
equal to the downwards closure of its {\it restriction} $${\mc
B}\upharpoonright M=\{t\in {\mc B}: t\subseteq M\}.$$ A further
examination of Rosenthal's result shows that for every pre-compact
family $\mc F$ of finite subsets of $\N$ there is an infinite set
$M$ such that the trace ${\mc F}[M]$ is actually equal to the
downwards closure of a uniform barrier $\mc B$ on $M$, or in other
words that the $\subseteq$-maximal elements of ${\mc F}[M]$ form a
uniform barrier on $M.$ As it turns out, this fact holds
considerably more information that the conclusion that ${\mc
F}[M]$ is merely a hereditary family which is especially
noticeable if one need to perform further refinements of $M$ while
keeping truck on the original family $\mc F.$ This observation was
the motivating point for our research which helped us to realize
that further extensions of Rosenthal's result require analysis of 
not only pre-compact families of finite subsets of $\N$ but also
maps from barriers into pre-compact families of finite subsets of
$\N,$ or, more generally, into weakly  compact subsets of $c_0$.
We have explained this point in our previous paper \cite{lop-tod},
where we have presented various results on partial 
unconditionality  such as near-unconditionality or
convex-unconditionality as  consequences of the structure theory
of this kind of mappings.  This paper is as a continuation of this
line of research. In  Section 3 we show how the combinatorics on
barriers can be used to prove the $c_0$-saturation for Banach
spaces $C(K)$ when $K$ is a countable compactum.  Recall that the
$c_0$-saturation of Banach spaces $C(K)$ over countable compacta
$K$ is a result originally due to Pe\l czy\'{n}ski and Semadeni
\cite{pel-sem} (see also \cite{arg-ka1} and \cite{ga-od-wa} for
recent accounts on this result.) More particularly, we show that
if $(x_i)\con C(K)$ is a normalized weakly-null sequence, then
there is $C\ge 1$, some infinite set $M$, some uniform barrier
$\mc B$ on $M$ of rank at most the Cantor-Bendixson rank of $K$
and some \emph{uniform} assignment $ \mu: \mc B\to c_{00}^+$ with
the property that $\supp \mu(s)\con s$ for every $s\in \mc B$, and
such that for every block sequence $(s_n)$ of elements of $\mc B$,
the corresponding sequence $(x(s_n))$ of linear combinations,
$$x(s_n)=\sum_{i\in s_n}(\mu(s_n))(i)x_i,$$ is a normalized block
sequence $C$-equivalent to the standard basis of $c_0$.

The last section  concerns the following natural measurement of unconditionality present in a given weakly
null sequence $(x_n)$ in a general Banach space $E.$ Given a family $\mc F$ of finite sets, we say that
$(x_n)$ is \emph{$\mc F$-unconditional} with constant at most $C\ge 1$ iff for every sequence of scalars
$(a_n)$,
$$\sup_{s\in \mc F}\nrm{\sum_{n\in s}a_n x_n}\le C \nrm{\sum_{n\in \N} a_n x_n}.$$
Thus, if for some infinite subset $M$ of $\N$ the trace $\mc F[M]$ contains the family of all finite subsets
of $M,$ the corresponding subsequence $(x_n)_{n\in M}$ is unconditional. Typically, one will not be able to
find such a trace, so one is naturally led to study this notion when the family $\mc F$ is pre-compact, or
equivalently, when $\mc F$ is a barrier. Since for every pair ${\mc F}_0$ and ${\mc F}_1$ of barriers on $\N$
there is an infinite set $M$ such that ${\mc F}_0[M]\subseteq{\mc F}_1[M]$ or ${\mc F}_1[M]\subseteq{\mc
F}_0[M]$ and since the two alternatives depend on the ranks of ${\mc F}_0$ and ${\mc F}_1,$ one is also
naturally led to the following measurement of unconditionality that refers only to a countable ordinal
$\gamma$ rather than a particular barrier of rank $\gamma.$ Thus, we say that a normalized basic sequence
$(x_n)$ of a Banach space $X$ is \emph{$\gamma$-unconditionally saturated} with constant at most $C\ge 1$ if
there is an $\gamma$-uniform barrier $\mc B$ on $\N$ such that for every infinite $M\con \N$ there is
infinite $N\con M$ such that the corresponding subsequence $(x_n)_{n\in N}$ of $(x_n)$ is $\overline{\mc
B\rest N}$-unconditional with constant at most $C$. (Here, $\overline{\mc B\rest N}$ denotes the topological
closure of the restriction ${\mc B\rest N}$ which in turn is equal to the trace $\mc B[N]$, a pleasant
property of any barrier.) It turns out that only indecomposable countable ordinals $\gamma$ matter for this
notion. We shall see, extending  the well-known example of Maurey-Rosenthal of a normalized weakly-null
sequence without unconditional subsequences, that every normalized basic sequence has a subsequence which is
$\omega$-unconditionally saturated, and that this cannot be extended further. For example, we show that for
every indecomposable countable ordinal $\gamma>\omega$ there is a compactum $K$ of Cantor-Bendixson rank
$\gamma +1$ and a normalized 1-basic weakly-null sequence $(x_n)\subseteq C(K)$ such that $(x_n)$ is
$\beta$-unconditionally saturated for all $\beta<\gamma$ but not $\gamma$-unconditionally saturated. More
precisely,  the summing basis of $c_0$ is finitely block-representable in every subsequence of $(x_n)$, and
so in particular, no subsequence of $(x_n)$ is unconditional.

\section{Preliminaries}
Let $\N$ denote the set of all non-negative integers and let
$\fin$ denote the family of all finite sets of $\N$. The topology
on $\fin$ is the one induced from the Cantor cube ${^\N}2$ via the
identification of subsets of $\N$ with their characteristics
function. Observe that this topology coincides with the one
induced by $c_0$, the Banach space of sequences converging to zero, with the same identification of finite sets and
corresponding characteristic functions.  Thus, we say that a family $\mc F\subseteq \fin$ is {\it compact} if it is a compact
space under the induced topology. We say that  $\mc F\subseteq
\fin$ is {\it pre-compact} if its topological closure
$\overline{\mc F}^{\mathrm{top}}$ taken in the Cantor cube
${^\N}2$ consists only of finite subsets of $\N.$ Given $X,Y \con
\N$ we write

\noindent (1) $X< Y$ iff $\max X<\min Y$. We will use the convention $\buit<X$ and $X<\buit$ for
every $X$.

\noindent (2)  $X \sqsubseteq Y$ iff $X \con Y$ and $X <
Y\setminus X$.

A sequence  $(s_i)$ of finite sets of integers is called a \emph{block sequence} iff $s_i<s_j$ for every
$i<j$, and it is called a \emph{$\De$-sequence} iff there is some finite set $s$ such that $s\sqsubseteq s_i$
($i\in \N$) and $(s_i\setminus s )$ is a block sequence. The set $s$ is called the \emph{root} of $(s_i)$.
Note that $s_i\to_i s$ iff for every subsequence of $(s_i)$ has a $\De$-subsequence with root $s$. It follows
that the topological closure $\overline{\mc F}$ of a pre-compact family $\mc F $ of finite subsets of $\N$ is
included in its downwards closure
$$\overline{\mc F}^{\con}= \conj{s\con t}{t\in \mc F}$$
with respect to the inclusion relation and also included in its downwards closure
$$\overline{\mc F}^\sqsubseteq= \conj{s\sqsubseteq t}{t\in \mc F}$$
with respect to the relation $\sqsubseteq.$ We say that a family $\mc F\con \fin$ is
$\subseteq$-{\it hereditary} if $\mc F=\overline{\mc F}^{\con}$ and $\sqsubseteq$-{\it hereditary}
if $\mc F=\overline{\mc F}^\sqsubseteq.$ The $\subseteq$-hereditary families will simply be called
{\it hereditary families}. We shall consider the following two restrictions of a given family $\mc
F$ of subsets of $\N$ to a finite or infinite subset $X$ of $\N$
\begin{align*}
\mc F\rest X = & \conj{s\in \mc F}{s\con X}, \\
\mc F [X]=& \conj{s\cap X}{s\in \mc F}.
\end{align*}

%
%
%
%
%
%
%
%

There are various ways to associate an ordinal index to a
pre-compact family $\mc F$ of finite subsets of $\N$. All these
ordinal indices are based on the fact that for $n\in \N$, the
index of the family
$$\mc F_{\{n\}}=
\conj{s\in \fin}{n<s, \, \{n\}\cup s\in \mc F}$$ is smaller or equal from that of $\mc F.$ For example, one
may consider the {\it Cantor-Bendixson index} $r(\mc F)$, the minimal ordinal $\alpha$ for which the iterated
Cantor-Bendixson derivative $\partial^{\alpha}(\mc F)$ is equal to $\emptyset$, then clearly $r(\mc
F_{\{n\}})\le r(\mc F)$ for all $n\in \N.$ Recall that $\partial\mc F$ is the set of all proper accumulation
points of $\mc F$ and that $\partial^{\alpha}(\mc F)=\bigcap_{\xi<\alpha}\partial(\partial^{\xi}(\mc F)).$
The rank is well defined since $\overline{\mc F}$ is countable and therefore a scattered compactum so the
sequence $\partial^{\xi}(\mc F)$ of iterated derivatives must vanish. Observe that if $\mc F$ is a nonempty
compact, then necessarily $r(\mc F)$ is a successor ordinal.

We are now ready to introduce the basic combinatorial concepts of
this section. For this we need the following piece of notation,
where $X$ and $Y$ are subsets of $\N$
$${_*}X=X\setminus \{\min X\}~ \mbox{and}~X/Y=\conj{m\in X}{\max
Y<m}$$ The set ${_*}X$ is called the \emph{shift} of $X$. Given
integer $n\in \N$, we write $X/n$ to denote $X/\{n\}=\conj{m\in
X}{m>n}$. The following notions have been introduced by
Nash-Williams.

\defi(\cite{nashwill})
Let  $\mathcal{F}\con \fin$.

\noindent (1)  $\mc F$ is called \emph{thin} if $s \not
\sqsubseteq t$ for every pair s, t of distinct members of
$\mathcal{F}$.

\noindent (2) $\mc F$ is called  \emph{Sperner} if $s \nsubseteq
t$ for every pair $s \neq t \in \mathcal{F}$.

\noindent (3) $\mathcal{F}$   is  called \emph{Ramsey}  if for
every finite partition $\mathcal{F} = \mathcal{F}_{0} \cup \cdots \cup \mathcal{F}_{k}
$ there is an infinite set $M \con \N$ such that at most one of the restrictions $\mc F_i \rest M$ is
non-empty.

\noindent (4)  $\mc F$ is called  a \emph{front} on $M$ if  $\mc F \con \mc P (M)$, it is thin, and for every infinite $N\con M$ there is
some $s\in \mc F $ such that $s \sqsubseteq N$.

\noindent (5)  $\mc F$ is called  a \emph{barrier} on $M$ if  $\mc F \con \mc P (M)$, it is Sperner, and for every infinite $N\con M$ there
is some $s\in \mc F $ such that $s \sqsubseteq N$.

\fdefi

Clearly, every barrier is a front but not vice-versa. For example, the family ${\N}^{[k]}$ of all $k$-element
subsets of $\N$ is a barrier. The basic result of Nash-Williams \cite{nashwill} says that every front (and
therefore every barrier) is Ramsey. Since as we will see soon there are many more barriers than those of the
form ${\N}^{[k]}$ this is a far reaching generalization of the classical result of Ramsey. To see a typical
application, let $\mc F$ be a front on some infinite set $M$ and consider its partition $\mc F = {\mc F}_0
\cup {\mc F}_1,$ where ${\mc F}_0$ is the family of all $\subseteq$-minimal elements of $\mc F$. Since $\mc
F$ is Ramsey there is an infinite $N\subseteq M$ such that one of the restrictions $\mc F_i \rest M$ is
empty. Note that $\mc F_1 \rest N$ must be empty. Since $\mc F_0 \rest N$ is clearly a Sperner family, it is
a barrier on $N$. Thus we have shown that every front has a restriction that is a barrier. Since barrier are
more pleasant to work with one might wonder why introducing the notion of front at all. The reason is that
inductive constructions lead more naturally to fronts rather than barriers. To get an idea about this, it is
instructive to consider the following notion introduced by Pudlak and R\"{o}dl.

\defi(\cite{pud-rod})
For a given countable ordinal $\al$, a family $\mc F$ of finite subsets of a given infinite set $M$ is called\emph{ $\al$-uniform on $M$} provided that:

\noindent(a) $\alpha =0$ implies $\mc F =\{\0\}$,

\noindent(b) $\alpha=\beta+1$ implies that  $\mc F_{\{n\}}$ is $\beta$-uniform on $M/n$,

\noindent(c) $\alpha > 0$ limit implies that there is an increasing sequence $\{\alpha_{n}\}_{n \in
M}$ of ordinals converging to $\alpha$ such that $\mc F_{\{n\}}$ is $\alpha_{n}$-uniform on $M/n$
for all $n \in M$.

$\mc F$ is called uniform on $M$ if it is $\al$-uniform on $M$ for some countable ordinal $\al$.
\fdefi

\nota\label{etirjhigs}

\noindent (a) If $\mc F$ is a front on $M$, then $\overline{\mc F}=\overline{\mc F}^{\ip}$.

\noindent (b) If $\mc F$ is uniform on $M$, then it is a front (though not necessarily a barrier)
on $M$.

\noindent (c) If $\mc F$ is $\al$-uniform (front, barrier) on $M$ and $\Theta:M\to N$ is the unique
order-preserving onto mapping between $M$ and $N$, then $\Theta"\mc F=\conj{\Theta"s}{s\in \mc F}$
is $\al$-uniform (front, barrier) on $M$.

\noindent (d) If $\mc F$ is $\al$-uniform (front, barrier) on $M$ then $\mc F\rest N$ is
$\al$-uniform (front, barrier) on $N$ for every $N\con M$.

\noindent (e) If $\mc F$ is  uniform (front, barrier) on $M$, then for every $s\in \overline{\mc
F}^{\ip}$ the family
$$\mc F_s=\conj{t}{s<t\text{ and } s\cup t\in \mc F}$$
is  uniform (front, barrier) on $M/s$.

 \noindent (e) If $\mc F$ is $\al$-uniform on $M$, then $\partial^{\al}(\overline{\mc
F})=\{\buit\}$, hence $r(\mc F)=\al+1$. (Hint: use that $\partial^{\be}(\mc
F_{\{n\}})=(\partial^{\be}(\mc F))_{\{n\}}$ for every $\be$ and every compact family $\mc F$).

\noindent (f) It is easy to prove by induction on $n$ that every $n$-uniform family on $M$ is of the form $M^{[n]}$. This is not the case in general.

\noindent (g)  An important example of a $\omega$-uniform barrier on $\N$ is the family $\mc S=\{s:
|s|=\min(s)+1\}$. We call $\mc S$ a {\it Schreier barrier} since its downwards closure is commonly called a
{\it Schreier family}. Indeed, it can be proved a $\mc B$ is a $\om$-uniform family on $M$ iff there is an unbounded mapping $f:M\to \om$ such that $\mc B=\conj{s\con M}{|s|=f(\min s)+1}$.
\fnota
The following result based on Nash-Williams' extension of Ramsey's theorem explains the relationship between
the concepts introduced above (see \cite{tod1} for proofs and fuller discussion).

\prop \label{teo5e} The following are equivalent for a family $\mathcal{F}$ of finite
subsets of $\N$:

\noindent (a) $\mathcal{F}$ is Ramsey.

\noindent (b) There is an infinite $M \con \N$ such that $\mathcal{F} \upharpoonright M$ is Sperner.

\noindent (c) There is an infinite $M \con \N$ such that $\mathcal{F} \upharpoonright M$ is either empty or uniform on $M$.

\noindent (d) There is an infinite $M \con \N$ such that $\mathcal{F} \upharpoonright M$ is either empty or a front on $M$.

\noindent (e) There is an infinite $M \con \N$ such that $\mathcal{F} \upharpoonright M$ is either empty or a barrier on $M$.

 \noindent (f)
There is an infinite $M \con \N$ such that $\mathcal{F} \upharpoonright M$ is thin.

\noindent (g) There is an infinite $M \con \N$ such that for every infinite $N \con M$ the restriction $\mathcal{F} \upharpoonright N$
cannot be split into two disjoint families that are uniform on $N$. \qed \fprop

In this kind of Ramsey theory one frequently performs  diagonalisation arguments  that can be
formalized using the following notion.

\defi\label{fusionsequence}
An infinite sequence $(M_k)_{k\in \N}$ of infinite subsets of $\N$ is called a \emph{fusion sequence of
subsets of $M\con \N$} if for all $k\in \N$:

\noindent (a) $M_{k+1}\con M_k\con M$,

\noindent (b) $m_k<m_{k+1}$, where $m_k=\min M_k$.

The infinite set $M_{\infty}=\{m_k\}_{k\in \N}$ is called the  \emph{fusion set (or limit)} of the sequence
$(M_k)_{k\in \N}.$
\fdefi

We have also the following simple facts connecting these combinatorial notions with the topological
concepts considered at the beginning of this section.

\prop \label{jjjfiiif} Fix a family $\mc F\con \fin$.

\noindent (a) If $\mc F$ is a barrier on $M$ then $\overline{\mc F}^\con=\overline{\mc
F}^\sqsubseteq=\overline{\mc F}$, and hence $\overline{\mc F}^\con$ is a compact family.


\noindent (b) If $\mc F$ is a barrier on $M$ then for every $N\con M$, $\overline{\mc F\rest
N}^\con=\overline{\mc F}^\con\rest N$.

\noindent (c) Suppose that $\mc F$ is a barrier on $M$. Then  for  every $N\con M$ such that $M\setminus N$
is infinite we have that $\mc F[N]=\overline{\mc F\rest N}^\con$, and in particular  $ \mc F[N]$ is downwards
closed.

\noindent (d) A family $\mc F\con M^{[<\infty]}$ is the
topological closure of a barrier on $M$ iff $\mc F^{\ip-\max}=\mc
F^{\con-\max}$ is a barrier on $M$.

\fprop

Barriers describe small families of finite sets, as it is shown in the following.

\teor\cite{lop-tod}\label{oherthhgrf}
Let $\mc F\con \fin$ be an arbitrary family. Then there is an infinite set $M\con \N$ such that either

\noindent (a) $\mc F[M]$ is the closure of a uniform barrier on $M$, or

\noindent (b) $M^{[\infty]}\con  \overline{\mc F}^\con$.
\fteor
Note that it follows that if $\mc F$ is pre-compact then condition (a) must hold.

We shall follow standard terminology and notation when dealing with sequences in Banach spaces (see
\cite{li-tza}). We recall now few standard definitions we are going to use along this paper.

\defi\label{mklsdjfoiews}
Let  $(x_i)$ be a sequence in a Banach space $E$.

\noindent (a) $(x_i)$ is called \emph{weakly-null} iff for every
$x^*\in E^*$, the sequence of scalars $(x^*(x_i))_i$ tends to 0.

\noindent (b) $(x_i)$ is called a \emph{Schauder basis} of $E$ iff  for every $x\in E$ there is a
unique sequence of scalars $(a_i)$ such that $x=\sum_i a_ix_i$. This is equivalent to say that
$x_i\neq 0$ for every $i$, the closed linear span of $(x_i)$ is $X$, and there is a constant
$\theta\ge 1$ such that for every sequence of scalars $(a_i)$, and every interval $I\con \N$,
\begin{equation}
\nrm{\sum_{i\in I} a_i x_i}\le \theta\nrm{\sum_{i\in \N} a_i x_i}.
\end{equation}

\noindent (c) $(x_i)$ is called a basic sequence iff it is a Schauder basis of its closed linear
span, i.e., $x_i\neq 0$ for every $i$, and there is $\theta\ge 1$ such that for every sequence of
scalars $(a_i)$, and every interval $I\con \N$, $\nrm{\sum_{i\in I} a_i x_i}\le \theta\nrm{\sum a_i
x_i}$. The infimum of those constants $\theta$ is called the basic constant of $(x_i)$.

\noindent (d) $(x_i)$ is called \emph{$\theta$-unconditional} ($\theta\ge 1$) iff for every
sequence of scalars $(a_i)$, and every subset $A\con \N$,
\begin{equation}
\nrm{\sum_{i\in A}a_ix_i}\le \theta\nrm{\sum_{i\in \N}a_ix_i}.
\end{equation}
$(x_i)$ is called unconditional if it is $\theta$-unconditional for some $\theta\ge 1$.

Given two basic sequences $(x_i)_{i\in M} $ and $(y_i)_{i\in N} $ of some Banach spaces $E$ and $F$, indexed
by the infinite sets $M,N\con \N$,  we say that  $(x_i)_{i\in M}\con E$ and $(y_i)_{i\in N}\con F$ are
$\theta$-equivalent, denoted by $(x_i)_{i\in M} \sim_\theta (y_i)_{i\in N}$,  if the order preserving
bijection $\Phi$ between the two index-sets $M$ and $N$ lifts naturally to an isomorphism between the
corresponding closed linear spans of these sequences sending $x_i$ to $ y_{\Phi(i)}$.
\fdefi

The \emph{ sequence of evaluation functionals} of $c_0$ is the biorthogonal sequence $(p_i)$ of the natural
basis $(e_i)$ of $c_0$, i.e. if $x=\sum_i a_ie_i\in c_0$, then $p_i(x)=a_i$. Note that weakly compact subsets
$K$ of $c_0$ are characterized by the property that every sequence in $K$ has a pointwise converging
subsequence to an element of $K$. It is clear that for every weakly-compact subset $K\con c_0$ the
restrictions of evaluation mappings $(p_i)$ to $K$ is weakly-null in $C(K)$. The sequence of restrictions
will also be denoted by $(p_i)$. Observe that $(p_i)$ as a sequence in the Banach space $C(K)$ is a monotone
basic sequence iff $K$ is closed under restriction to initial intervals.

There are two particularly important examples of weakly-compact subsets of $c_0$ naturally associated to a  normalized weakly null sequence $(x_i)_{i\in M}$ of a Banach space $E$:

\noindent (a) the set
$$R_E((x_i)_{i\in M})=\conj{(x^*(x_i))_{i\in M}\in c_0}{x^*\in B_{E^*}}$$
 is  symmetric, 1-bounded and  weakly-compact subset of
$c_0$.

\noindent (b) If   $E=C(K)$, $K$ compactum, then the set
$$R_K((x_i)_{i\in M})=\conj{(x_i(c))_{i\in M}\in c_0}{c\in K}$$
is also  1-bounded and weakly-compact.

In both cases one has that  $(x_i)_{i\in M}$ is 1-equivalent to the evaluation mapping sequences of $C(R_E((x_i)_{i\in M}))$ and $C(R_K((x_i)_{i\in M}))$.

We say that a subset $X$ of $c_0$ is {\it weakly pre-compact} if its closure relative to the weak topology of
$c_0$ is weakly compact.  We have then the following, not difficult to prove.
\prop

\noindent (a) $\mc F\con \fin $ is pre-compact iff the set $\conj{\chi_s}{s\in \fin}\con c_0$ of
characteristic functions of sets in $\mc F$ is weakly-pre-compact.

\noindent (b) For every weakly-pre-compact subset $X$ of $c_0$ and every $\vep>0$ one has that
$$\supp_\vep X=\conj{\conj{n\in \N}{|\xi(n)|\ge \vep}}{\xi \in X} \text{ is pre-compact}.$$

\fprop

Finally, we introduce few combinatorial notions concerning mappings from families of finite sets of integers
into $c_0$. For more details see \cite{lop-tod}.
\defi (\cite{lop-tod})
Let $\mc F\con \fin$ be an arbitrary family, and let $f:\mc F\to c_0$.

\noindent (a) $f$ is internal if for every $s\in \mc F$ one has that $\supp f(s)\con s$.

\noindent (b) $f$ is uniform if for every $t\in \fin$ one has that
$$|\conj{\vphi(s)(\min (s/t))}{t\ip s , \, s\in \mc F}|=1$$

\noindent (c) $f$ is Lipschitz if for every $t\in \fin$ one has that
$$|\conj{\vphi(s)\rest t }{t\ip s , \, s\in \mc F}|=1$$

\noindent (d) $f$ is called a $U$-mapping if $\mc F$ if it is internal and uniform.

\noindent (e) $f$ is called a $L$-mapping if $\mc F$ if it is internal and  Lipschitz.
\fdefi

\nota
\noindent (a) Every uniform mapping is Lipschitz, but the reciprocal  is in general false. For example, the
mapping $f:\fin\to c_0$ defined by $f(s)(i)=i$ if $i\in s$ and $f(s)(i)=0$ is Lipschitz but not uniform.

\noindent (b) Every $L$-mapping $f:\mc F\to c_0$ can be naturally extended to a continuous mapping
$f':\overline{\mc F}^{\ip}\to c_0$ by setting $f'(t)=f(s)\rest t$ for (any) $s\in \mc F$   such that $t\ip
s$.

 \noindent (c)  The importance of internal mappings can be seen, for example, by the well-known
result of Pudlak-R\"{o}dl \cite{pud-rod} stating that if $f:\mc B\to X$ is a function defined on a barrier $\mc
B$ on $M$ then there is $N\con M$, a barrier $\mc C$ on $N$, and an internal mapping $g:\mc B\rest N\to \mc
C$ such that for every $s,t\in \mc B\rest N$ one has that $f(s)=f(t)$ iff $g(s)=g(t)$.

\noindent (d)   $U$-mappings were used in \cite{lop-tod} to produce some weakly-null sequences playing
important role in the better understanding of an abstract concept of unconditionality (see \cite{lop-tod} for
more details).
\fnota
The main result on mappings
defined on barriers is the following:
\teor\cite{lop-tod}\label{jhhuheurx}
Suppose that $\mc B$ is a barrier on $M$, $K\con c_0$ is weakly-compact and suppose that $f:\mc B\to K$. Then
for every $\vep>0$  there is $N\con M$  and there is a $U$-mapping $g:\mc B\rest N\to c_{00}$ such that for
every $s\in \mc B\rest N$ one has that
$$\nrm{f(s)\rest N -g(s)}_{\ell_1}\le \vep.$$
\fteor
\cor\label{jsrhjrhfd}
Suppose that $f:\mc B\to c_0$ is an internal mapping defined on a barrier $\mc B$. Suppose that in addition
$f$ is bounded, i.e. there is $C$ such that for every $s\in \mc B$ one has that $\nrm{f(s)}_\infty \le C$.
Then for every $\vep>0$ there exists is a $U$-mapping $g:\mc B\rest N\to c_{00}$ such that for every $s\in
\mc B\rest N$ one has that
$$\nrm{f(s)-g(s)}_{\ell_1}\le \vep.$$
\fcor
\prue
Let us prove first that the image of $f$ is weakly-pre-compact: For suppose that $(f(s_n))_n$ is an arbitrary
sequence. Let $M\con \N$ be such that $(\supp f(s_n))_{n\in M}$ converges to some $s\in \overline{\mc
B}^\ip$. This is possible because $f$ is internal. Since $f$ is bounded, we can find $N\con M$ such that
$(f(s_n))_{n\in N}$ is weak-convergent in $c_0$.

Now the desired result follows from \ref{jhhuheurx} by using that $f$ is in addition internal.
\fprue

\section{$c_0$-saturation of $C(K)$ for a countable compactum $K$}

Recall the result of Pelczynski and Semadeni \cite{pel-sem}     which says that every Banach space of the
form $C(K)$ for $K$ a countable compactum is $c_0$-saturated in the sense that every of its closed
infinite-dimensional subspaces contains an isomorphic copy of $c_0.$ The purpose of this section is to
examine the $c_0$-saturation using the theory of mappings on barriers developed above in Section 3. We start
with a convenient reformulation of the problem. We start with a definition.
\defi For a given subset $X$ of $c_0$, let $\supp X=\conj{\conj{i\in \N}{\xi(i)\neq 0}}{\xi\in X}$ be the \emph{support set} of $X$.
We say that a weakly compact subset $K$ of $c_0$ is \emph{supported by a barrier on $M$} if  its support set $\supp K$ is the
is the closure of a uniform barrier on $M$.
\fdefi

\lema\label{jjkweuiorghh}
Suppose that $K$ is a countable compactum. Suppose that $(x_i)\con C(K)$ is a normalized weakly
null sequence. Then for every $\vep>0$ there is  subsequence $(x_i)_{i\in M}$ and a weakly-compact
subset $L\con c_{0}$  supported by a barrier on $\N$ of rank not bigger than the Cantor-Bendixson rank of $K$ such that
 $(x_i)_{i\in M}$ and the evaluation mapping $(p_i)_{i\in \N}$ of $C(L)$ are $(1+\vep)$-equivalent.
\flema
\prue
Fix $\vep>0$.   Find first an strictly decreasing sequence $(\vep_i)$ such that $\sum_i \vep_i\le \vep$ and such that
\begin{equation}\label{ouetoioidsdjfs}
\conj{\vep_i}{i \in \N}\cap \conj{|x_i(c)|}{c\in K}=\buit.
\end{equation}
This is possible because $K$ is countable. Now define $\vphi:K\to \mc P(\N)$ by $\vphi(c)=\conj{i\in \N}{|x_i(c)|\ge \vep_i}$. Note that \eqref{ouetoioidsdjfs} implies that $\vphi$ is a continuous function. Enumerate $K=\{c_k\}_{k\in \N}$. Since $(x_i)$ is weakly-null we can find  a fusion sequence $(M_k)$ such that for every $k$  and    every $i\in M_k$ one has that $|x_i(c_k)|<\vep_k$. Now if we set $N$ to be the corresponding fusion set  then for every $k$ one has that $\conj{i\in M}{|x_i(c_k)|\ge \vep_i}\con \{n_0,\dots,n_{k-1}\}$. This means that the mapping $\psi=\xi_M \cdot \psi$ is continuous with image included in $\fin$.  Set $N={_*}M$ and denote the immediate predecessor  of $i\in N$ in $M$ by $i^-$.
Since $K$ is a zero-dimensional compactum, we can find clopen sets $C_i\con K$ $(i\in N)$ such that
$$\text{$K\setminus x_i^{-1} ((-\vep_{i^-} ,\vep_{i^-}))\con   C_i \con K\setminus  x_i^{-1} ([-\vep_{i} ,\vep_{i}]) $ for every
$i\in N$. }$$
Set $y_i=\chi_{C_i}x_i$ for each $i\in N$.  So one has

\noindent(i) $\nrm{x_i-y_i}_{K}<\vep_{i^-}$, so  $(x_i)_{i\in N}$ and $(y_i)_{i\in N}$ are
$1+\vep$-equivalent, and

\noindent(ii) for every $c\in K$ and every $i\in N$, if $|
y_i(c)|\le \vep_i$, then $y_i(c)=0$.

Since for every $c\in K$, by (ii) above, one has that
$$\conj{i\in N}{y_i(c)\neq 0}=\conj{i\in N}{c\in C_i \text{ and } |x_i(c)|\ge \vep_i}=\psi(c),$$
it follows that the support set  $\mc F$ of $R_K((y_i)_{i\in N})$   coincide with the image of $\psi$, so it is a compact family of $\N$.  We use now Theorem \ref{oherthhgrf} to find  $P\con N$ such that $\mc F[P]$ is the closure of a uniform barrier on $P$. This implies that $R_{K}((y_i)_{i\in P})$ is supported by  a barrier $\mc B$ on $P$.    Let $\theta$ be the unique order preserving mapping from $\N$ onto $P$, and let  $ \Theta:  c_0\rest P=\conj{\xi\in c_0}{\supp \xi\con P}\con c_0 \to c_0$ be defined by $\Theta(\xi)(n)=\xi(\theta(n))$.
 This is an homeomorphism between $c_0\rest P$ and $c_0$, both with the weak topology, so    $L=\Theta" R_{K}((y_i)_{i\in P})$  is a weakly-compact subset of $c_0$ and supported by the barrier $\theta^{-1}\mc B=\conj{\theta^{-1}s}{s\in \mc B}$ on $\N$.  Now it is easy to see that the evaluation mapping $(p_i)_{i\in \N}$ of $C(L)$ is a normalized weakly-null sequence $1+\vep$-equivalent to $(x_i)_{i\in P}$.
\fprue

\teor
 Suppose that
$(x_i)\con C(K)$ is a normalized weakly-null sequence for a  countable compactum  $K$. Then there is a constant $C\ge 1$,  an infinite set $M$, a uniform
barrier $\mc B$ on $M$  whose rank  is at most the Cantor-Bendixson rank of $K$, and some $U$-mapping  $ \mu: \mc
B\to c_{00}^+$ such that  for every block sequence $(s_n)\con \mc B$ the corresponding sequence of linear combinations $(\sum_{i\in s_n}(\mu(s_n))(i)x_i)_n$ is a
normalized block sequence $C$-equivalent to the unit vector basis of $c_0$.

\fteor
\prue The proof is by induction on the Cantor-Bendixson rank of $K$. First of all, by Lemma
\ref{jjkweuiorghh} we may assume that $K$ is a weakly-compact subset of $c_0$ supported by a barrier $\mc B$ on $\N$ and that  the normalized weakly null sequence $(x_i)$ is the corresponding evaluation mapping sequence $(p_i)_{i\in \N}$. If
$\al=1$, then $\mc B=\N^{[1]}$ and clearly $(p_i)$ is equivalent to the unit vector basis of $c_0$.
So assume that $\al>1$. By going to a subsequence of $(p_i)$ if needed, we may also assume in this
case that $|s|\ge 2 $ for every $s\in \mc B$.  For each integer $n$ set $\mc F_n=\overline{\bigcup_{m\le n}\mc
B_{\{m\}}}$. Since $\mc B$ is a $\al$-uniform family, we have that for every $n$,
$\partial^{\al}\mc F_n=\buit$, so its Cantor-Bendixson rank is strictly smaller than $\al+1$.
For each $n\in \N$, let
$$K_n=\conj{f\rest s}{s\in \mc F_n}.$$
This is a compactum   whose support is $\mc F_n$  and whose rank
 is strictly smaller than $\al+1$. So, the evaluation mapping sequence $(p_i)$ is a weakly-null
 sequence of $C(K_n)$ for every $n$. Observe that for every sequence of scalars $(a_i)$ we have
 that
\begin{equation}
 \nrm{\sum_i a_ip_i}_{n}= \nrm{\sum_i a_ip_i}_{K_n}=\sup \conj{\nrm{\sum_{i\in s} a_i p_i}_{K}}{s\in \mc F_n}.
 \end{equation}
Using the fact that the family $\mc F_n$ is hereditary, we obtain that $(p_i)$ is 1-unconditional.  Since we assume
that all the singletons $\{i\}$  belong to $\mc F_{n}$, it follows  that $(p_i)_{i\ge 1}$ is  indeed a
1-unconditional normalized weakly null sequence in $C(K_n)$.

Fix $\vep>0$, and let $(\vep_n)_n$ be a summable sequence with $\sum_n \vep_n<\vep/2$.   By the Ramsey
property of the uniform barrier $\mc B$, we can find a fusion sequence $(M_k)_k$ such that, setting $n_k
=\min M_k $ for each $k\in \N$, we have that for every $k$ the following dichotomy holds:

(I) Either for every $s\in \mc B\rest M_{k}$ there is some $\mu_k(s)\in c_{00}$  with $\supp \mu_k(s)\con s$,
$0\le \mu_k(s)(i)\le 1$ for every $i\in s$, and such that for every such that $\nrm{\sum_{i\in s}\mu_k(s)(i)
p_i}_K=1$ while $\nrm{\sum_{i\in s}\mu_k(s)(i) p_i}_{n_k}<\vep_k$, or else

(II) $\nrm{\sum_{i\in s}a_i p_i}_K\le 2\vep_k^{-1} \nrm{\sum_{i\in s}a_i p_i}_{n_k}$ for every $s\in \mc
B\rest M_{k}$ and every $(a_i)_{i\in s}$.

Suppose first that (I) holds for every $k$.  Let $M_{\infty}=\{n_k\}$ be the corresponding fusion set. Then
let $\mc C=\mc B\rest M_\infty$. For $s\in \mc C$, define $\mu(s)=\mu_k(s)$, where $n_{k}=\min s $.
This is well defined since $s\in \mc B \rest M_{k}$.  For  a given $s\in \mc C$, let
$$x(s)=\sum_{i\in s}\mu(s)(i)p_i.$$
Our intention is to show that for every block sequence $(s_i)_i$ in $\mc C$ one has that   $(x(s_i))_{i}$ is
$2+\vep$-equivalent to the $c_0$-basis.     So fix such sequence $(s_i)$ and let  $(b_i)_{i\in \N}$ be a
sequence of scalars with $|b_i|\le 1$ for every integer $i$.   Since each $x(s_i)$ is normalized and since
$(p_i)$ is monotone, we obtain that
$$\nrm{\sum_i b_i x(s_i)}_{K}\ge (1/2)\nrm{\sum_i b_i e_i}_\infty.$$
Suppose that $\xi\in K$, and let $i_0=\min\conj{i }{ s_i \cap \supp \xi \neq \buit}$. Fix $i>i_0$,
and let $k_i$ be such that $n_{k_i}=\min s_i$. Since  $\supp \xi \cap s_i\in \mc F_{\max s_{i_0}}$
we have that
\begin{equation}
| x(s_i) (\xi)|\le \nrm{\sum_{j\in  s_i \cap \supp f} a_i^{(k_i)}p_i}_{\max s_{i_0}}<
\vep_{k_i}.
\end{equation}
It follows that
\begin{equation}
| \sum_i b_i x(s_i)(\xi)| \le |b_{i_0}| +\sum_{i>i_0}|b_i||x(s_i)(\xi)| \le |b_{i_0}|+
\frac\vep2.
\end{equation}
So, $\nrm{\sum_{i}b_i x(s_i)}_{K}\le (1+\vep/2)\nrm{\sum b_i e_i}_{\infty}$. Finally use Corollary
\ref{jsrhjrhfd} to perturb $\mu$ and make it $U$-mapping.

Suppose now that $k_0$ is the first $k$ such that (II) holds for $k$. Set $M=M_{k}$. It readily follows that
for every $x$ in the closed linear span of $(p_i) _{i\in M}$ one has that $\nrm{x}_{K}\le
\vep_{k_0}^{-1}\nrm{x}_{n_{k_0}}$.  By inductive hypothesis applied to $(p_i)\con C(K_{n_{k_0}})$, there is
some $C\ge 1$, some uniform barrier $\mc C$ on some $N\con M$ of rank not bigger than the one of
$K_{n_{k_0}}$ and some $\mu$ fulfilling the conclusions of the Lemma. Fix $s\in \mc C$. Then
$\nrm{\mu(s)}_{n_{k_0}}=1$, so we can find some $t_s\con s$ such that
$1=\nrm{\mu(s)}_{n_{k_0}}=\nrm{\mu(s)\rest t_s}_{K}$. Observe that, by 1-unconditionality of
$\nrm{\cdot}_{n_{k_0}}$,
 $\nrm{\mu(s)\rest t}_{ n_{k_0}}=1$. Define $\nu: \mc C\to c_{00}$ by $\nu(s)=\mu(s)\rest t_s$.   Finally, let us check that $(x(s_i))\con C(K)$ is
 $C\vep_{n_k}^{-1}$-equivalent
 to  the $c_0$-basis for every block sequence $(s_i)_i$ in $\mc C$.  Fix scalars $(a_i)$, $|a_i|\le 1$ ($i\in \N$).
   We obtain the inequality $\nrm{\sum_i a_i \nu(s_i)}_{K}\ge (1/2)
 \nrm{\sum_i a_i
 e_i}_\infty$  by the monotonicity of the basic sequence $(p_i)$. Now,
\begin{equation}
\nrm{\sum_i a_i \nu(s_i)}_{K}\le\frac1{ \vep_{n_{k_0}}}\nrm{\sum_i a_i
\nu(s_i)}_{n_{k_0}}\le \frac1{ \vep_{n_{k_0}}}\nrm{\sum_i a_i \mu(s_i)}_{n_{k_0}} \le
\frac{C}{ \vep_{n_{k_0}}}\nrm{\sum_i a_i e_i}_\infty.
\end{equation}
\fprue

\section{Conditionality}

We start with the following natural slightly variation on the notion of $S_\xi$-unconditionality from \cite{arg-god-ros}, and which is a generalization of unconditionality  (see Definition \ref{mklsdjfoiews} (d)).

\defi
 Let $\mc F$ be a family of finite sets of integers. A normalized basic sequence $(x_n)$ of a Banach space $E$ is called
\emph{$\mc F$-unconditional}
with constant at most $C\ge 1$ iff for every sequence of scalars $(a_n)$,
$$\sup_{s\in \mc F}\nrm{\sum_{n\in s}a_n x_n}\le C \nrm{\sum_{n\in \N} a_n x_n}.$$
\fdefi
This generalizes the notion of unconditionality covered by the case of $\mc F=\fin$. The question is whether
every normalized weakly-null sequence has a $\mc F$-unconditional subsequence. Observe that  the subsequence
$(x_n)_{n\in M}$ is $\mc F$-unconditional iff  it is $\mc F[M]$-unconditional, so the existence of an $\mc
F$-unconditional subsequence is closely related to the form of the traces $\mc F[M]$. If we assume that in
addition the family $\mc F$ is hereditary, then, by the Theorem \ref{oherthhgrf}, two possibilities can
occur: The first one is that some trace of $\mc F$ consists on all finite subsets of some infinite set $M$.
In this case, for subsequences of $(x_n)_{n\in M}$ the $\mc F$-unconditionality coincides with the
unconditionality. The second case is  when some trace of $\mc F$ is the closure of a uniform barrier. So one
is naturally led to examining the standard compact families of finite subsets of $\N$.   We begin
with the following positive result announced in \cite{mau-ros} and first proved by E. Odell \cite{odell1}
concerning the Schreier family $\overline{\mc S}=\{s\subseteq {\N}:|s|\leq \min(s)+1 \}$.
\teor
\label{schreier} Suppose that $(x_n)$ is a  normalized weakly-null sequence of a Banach space $E$. For every
$\vep>0$ there is  a $\overline{\mc S}$-unconditional subsequence with constant $2+\vep$. \qed
\fteor

%

Recall that if $\mc F$ is a barrier on some set $M$ then its trace ${\mc F}[N]$ on any co-infinite subset $N$
of $M$ is hereditary and that for every pair ${\mc F}_0$ and ${\mc F}_1$ of barriers on the same domain $M$
there is an infinite set $N\subseteq M$ such that ${\mc F}_0[N]\subseteq{\mc F}_1[N]$ or ${\mc
F}_1[N]\subseteq{\mc F}_0[N]$. Since the two alternatives are dependent on the ranks of ${\mc F}_0$ and ${\mc
F}_1,$ one is naturally led to the following measurement of unconditionality.

\defi
Suppose that $\ga$ is a countable ordinal. A normalized basic
sequence $(x_n)$ of a Banach space $E$ is called
\emph{$\ga$-unconditionally saturated} with constant at most $C\ge
1$ if for every $\ga$-uniform barrier $\mc B$ on $\N$ and for
every infinite $M$ there is infinite $N\con M$ such that the
corresponding subsequence $(x_n)_{n\in N}$ of $(x_n)$ is
$\overline{\mc B}$-unconditional with constant at most $C$.

We say that $(x_n)_n$ is $\ga$-unconditionally saturated if it is $\ga$-unconditionally saturated with constant $C$ for some $C\ge 1$.
\fdefi

\nota

\noindent (a) A sequence $(x_n)_n$ is $\ga$-unconditionally saturated iff  given a $\ga$-uniform barrier $\mc
B$ every subsequence of $(x_n)_n$ has a further $\mc B$-unconditional subsequence. The reason for this  is
that given two $\ga$-uniform barriers $\mc B$ and $\mc C$ on a set $M$ we have that there is $N\con M$ such
that either $\overline{\mc B\rest N}\con \overline{\mc C\rest N}\con \overline{\mc B\rest N}\oplus N^{[\le
1]}$ or the symmetric situation holds, where $\mc F\oplus \mc G=\conj{s\cup t}{s\in \mc G,\, t\in\mc F \text{
and }s<t}$ (see \cite{tod1}).

\noindent (b)  It follows from Theorem \ref{schreier} that every normalized weakly null sequence is
$\om$-unconditionally saturated.  Since the $\om$-uniform barriers are of the form $\conj{s\in
\fin}{|s|=f(\min s)+1}$ for some unbounded mapping $f:M\to \N$ one can easily modify the proof of Theorem
\ref{schreier} to prove that  every normalized weakly-null sequence is  $\om$-unconditionally saturated with
constant at most $2+\vep$.

\noindent (c) If the normalized basic sequence $(x_n)$ is monotone, then it is $\overline{\mc
B}$-unconditional iff it is $\mc B$-unconditional for every uniform barrier  $\mc B $ on $\N$.

\noindent (d) An analysis of the Maurey-Rosenthal \cite{mau-ros} example of a weakly-null sequence $(x_n)$
with no unconditional basic subsequence (see Example \ref{mvbjfkeweer} below) reveals an $\omega^2$-uniform
barrier $\mc B_{\mathrm{MR}}$ such that no infinite subsequence $(x_n)_{n\in M}$ is $\mc
B_{\mathrm{MR}}$-unconditional with any finite constant $C$. So this is an example of a normalized
weakly-null sequence with no $\om^2$-unconditionally saturated subsequence.

\noindent (e) Recall that an ordinal $\ga$ is called indecomposable if for every $\be<\ga$, $\be\om\le \ga$.
Equivalently, $\ga=\om^{\be}$ for some $\be$.   Suppose that $\ga$ is the maximal indecomposable ordinal
smaller than some fixed ordinal $\al$. Then a normalized basic sequence $(x_n)$ is $\al$-unconditionally
saturated if and only it is $\ga$-unconditionally saturated.
\fnota

 \eje\label{mvbjfkeweer} First of all, for a fixed
$0<\varepsilon<1$ choose a fast increasing sequence $(m_{i})$ such that
\begin{equation}
\sum_{i=0}^{\infty}\sum_{j\neq i}
\min\{(\dfrac{m_{i}}{m_{j}})^{1/2}, (\dfrac{m_{j}}{m_{i}})^{1/2}\}
\leq \frac\varepsilon{2}.
\end{equation}
Let  $ \fin ^{[<\infty]}$ be the collection of all finite block sequences
$\mathit{E_{0}<E_{1}<\dots<E_k}$ of nonempty finite subsets of $\N$. Now  choose a $1-1$ function
\begin{equation}
\sigma :  \fin ^{[<\infty]} \rightarrow \{m_i\}
\end{equation} such that $\varphi((s_i)_{i=0}^n)>s_n$ for all $(s_i)\in  \fin ^{[<\infty]}$
Now let $\mc B_{\mathrm{MR}}$ be the family of unions $s_0\cup s_1\cup \dots\cup s_n$ of finite sets such
that

\noindent (a) $(s_i)$ is block and  $s_0=\{n\}$.

\noindent (b) $|s_i|=\sig(s_0,\dots,s_{i-1})$  ($1\le i\le n$).

It turns out that $\mc B_{\mathrm{MR}}$ is a  $\omega^2$-uniform barrier on $\N$ (see Proposition
\ref{dfjklweobf} below), hence $\overline{\mc B_{\mathrm{MR}}}=\overline{\mc
B_{\mathrm{MR}}}^{\sqsubseteq}$ is a compact family with rank $\om^{2}+1$. Observe that by
definition, every $s\in \mc B_{\mathrm{MR}}$ has a unique decomposition $s=\{n\}\cup s_1\cup \dots
\cup s_n$  satisfying (a) and (b) above. Now define the mapping $\Phi:\mc B_{\mathrm{MR}}\to
c_{00}$,
\begin{equation}
\Phi(s)= e_{n} +\sum_{i=1}^n \frac{1}{|s_i|^{\frac12}}\sum_{k\in s_i}e_k.
\end{equation}
It follows that $\Phi$ is a $U$-mapping defined on the barrier $\mc B_{\textrm{MR}}$.   Now we can define the
Banach  space $\mk X_{\textrm{MR}}$ as the completion of $c_{00}$ under the norm
$$\nrm{x}_{\textrm{MR}}=\sup\conj{|\langle \Phi(s),x \rangle|}{s \in \mc B_{\textrm{MR}}}.$$
The natural Hamel basis $(e_n)$ of $c_{00}$ is now a normalized weakly-null monotone basis of $\mk
X_\textrm{MR}$ without unconditional subsequences. Indeed, without $\om^2$-unconditionally saturated
subsequences. Moreover this  weakly-null sequence has the property that the summing basis $(S_i)$ of $c$, the
Banach space of \emph{convergent sequences of reals}, is finitely-block representable in the linear span of
every subsequence of $(e_i)$ (and so the summing basis of $c_0$), more precisely, for every $M$, every $n\in
\N$ and every $\vep>0$ there is a normalized block subsequence $(x_i)_{i=0}^{n-1}$ of $(e_i)_{i\in M}$ such
that for every sequence of scalars $(a_i)_{i=0}^{n-1}$,
\begin{equation*}
\max\conj{|\sum_{i=0}^{m}a_i|}{m<n}\le \nrm{\sum_{i=0}^{n-1}a_ix_i}_{C(K)}\le
(1+\vep)\max\conj{|\sum_{i=0}^{m}a_i|}{m<n}.
\end{equation*}
On the other hand, by Proposition \ref{schreier} the sequence $(p_i)$ is $\om$-unconditionally saturated
with constant $\sim 2$.

Another presentation of this space is the following: Since $\Phi$ is uniform, it is Lipschitz, so there is a
unique extension $\Phi:\overline{\mc B_{\mathrm{MR}}}\to c_{00}$, naturally defined by $\Phi(s)=\Phi(t)\rest
t$, where $t\in \mc B_{\mathrm{MR}}$ is (any) such that $s \sqsubseteq t$. Now define $K=\Phi"\mc
B_{\mathrm{MR}}\con c_{00}$. This is a weakly-compact subset of $c_{00}$ whose  rank the same than
$\overline{\mc B_{\mathrm{MR}}}$, i.e., $\om^2+1$. Then  the corresponding evaluation sequence $(p_i)\con
C(K)$ is $1$-equivalent to the basis $(e_i)_i$ of $\mk X_{\textrm{MR}}$. \feje

Building on the idea of   Example \ref{mvbjfkeweer}, we are now going to find, for every countable indecomposable ordinal
$\ga$, a $U$-sequence   with no unconditional subsequences but
$\be$-unconditionally saturated for every $\be<\ga$.  Before embarking into the construction, we need to
recall a localized version of Pt\'{a}k's Lemma. For this we need the following notation: Given a family $\mc F$,
and $n\in \N$, let
$$\mc F \otimes n=\conj{s_0\cup \dots \cup s_{n-1}}{(s_i)_{i=0}^{n-1}\con \mc F\text{ is
block}}.$$
 It can be shown that $\mc F\otimes n$ is a $\al n$-uniform family if $\mc F$ is an
$\al$-uniform family.

Given $\xi\in c_{00}$ we will write $\xi^{1/2}$ to denote $(\xi(i)^{1/2})$. Given $\xi\in c_{00}$
and a finite set $s$, let $\langle \xi , s\rangle =\langle \xi,\chi_s\rangle=\sum_{i\in s}\xi(i)$.

\defi
A \emph{mean} is an element  $\mu\in c_{00}^+$  with the property that $\sum_{i\in \N}\mu(i)=1$.  We say that
$\mu: \mc B\to c_{00}^+$  is a \emph{$U$-mean-assignment} if $\mu$ is a $U$-mapping such that for every $s\in
\mc B$ one has that $\mu(s)$ is a mean.
\fdefi

\lema Suppose that $\mc B$ is an $\al$-uniform barrier on $M$,
$\al\ge 1$.  Let $\ga=\ga(\al)$ be the maximal indecomposable ordinal not bigger than $\al$,and let
$n=n(\al)\in \N$, $n\ge 1$, be such that $\ga n\le \al< \ga (n+1)$. Then for every $k\in \N$, $k>1$,  every
$\vep>0$, and every $\be$-uniform barrier $\mc C$ on $M$ with $\be> \al k$ there $N\con M$ and
$U$-mean-assignment $\mu:\mc C\rest N\to c_{00}^+$ such that

\begin{equation}\label{kghhthh}
\sup\conj{\langle\mu(s)^\frac12,t\rangle }{t\in \mc B}\le \frac{ (1+\vep)(n+1)}{(n k)^{\frac12}}
\end{equation}
for every $s\in \mc C\rest N$.
\flema
\prue
The proof   is by induction on $\al$. Fix $\vep>0$ and $k>1$.  Let $\mc C$ be an $\be$-uniform
family on $M$ such that $\be>\al k$.

Notice that  if we prove that for every $N\con M$ there is  one mean $\mu$ with support in $ \mc C\rest N$
such that \eqref{kghhthh} holds, then the   Ramsey property of the uniform barrier $\mc C$ gives the
existence for some $N\con M$ of a mean-assignment $\mu:\mc C\rest N\to c_{00}$ such that $\mu(s)$ has the property   \eqref{kghhthh}
for every $s\in \mc C \rest N$. Then Corollary \ref{jsrhjrhfd} gives the desired $U$-mapping.

Let $\mc D$ be a $\ga$-uniform barrier on $M$ (if $n=1$ we take $\mc D=\mc B$), and fix $N\con M$.
Find first $P\con N$ be such that $(\mc D\otimes nk) \rest P\con \overline{\mc C}$ as well as $\mc
B\rest P\con \overline{\mc D\otimes (n+1)}$. Consider $(\ga_i)_{i\in P}$ such that $\mc
D_{\{i\}}\rest P$ is $\ga_i$-uniform on $P/i$. Observe that for every $i\in P$  we have that $\ga_i
<\ga$, so, since $\ga$ is indecomposable, $\ga_i \om \le \ga$. Let $\mu_0$ be any mean such that
$\supp \mu_0\in \mc B\rest P$. By inductive hypothesis applied to appropriate $\al_i$'s, we can
find a block sequence $(\mu_j)_{j=0}^{nk-1}$ of means with support in $\mc B\rest P$ such that for
every $1\le j\le nk-1$,
\begin{equation}\label{vmnfhrer}
\sup\conj{\langle\mu_j^{\frac12}, t\rangle}{t\in \mc D, \text{ and $\min t\le \max \supp
\mu_{j-1}$}}<\frac{\vep}{2^{j+1}}.
\end{equation}
Let  $\nu=(1/(nk))\sum_{j=0}^{nk-1}\mu_j$. Observe that $\supp \nu\in \overline{(\mc D \otimes (nk))\rest
P}\con\overline{\mc C}$. Then, for every $t\in \mc B$, by \eqref{vmnfhrer},
\begin{equation}\label{vbftwww}
\langle\nu^{\frac12},t\rangle = \frac1{(nk)^{\frac12}}\sum_{j=0}^{k-1} \sum_{i\in t}\mu_j(i)^\frac12 \le
  \frac{1+\frac\vep2}{(nk)^{\frac12}}.
\end{equation}
Let us point out that $\supp \nu$ is, possibly, not a set in $\mc C$. However  it is easy to slightly perturb
$\nu$  to a newer mean with support in $\mc C$ and satisfying \eqref{vbftwww} for every $t\in \mc B$: Let
$s\in \mc C$ be such that $\supp \nu\ip s$, and set $u=s\setminus \supp \nu$.  Let $\de>0$ be such that
\begin{equation}
(1+\frac\vep2)(1-\de)^{1/2} +(nk\de |u|)^{1/2} \le 1+\vep.
\end{equation}
Now set
\begin{equation}
{\mu}= (1-\de)\nu+ \frac{\de}{|u|} \chi_{u}.
\end{equation}
${\mu}$ is a mean whose support is $s\in \mc C$. It can be shown now that for every $t\in \mc B$,
\begin{equation}
\sum_{i\in t} {\mu}(i)^{\frac12} \le      \frac{1+\vep}{(nk)^\frac12},
\end{equation}
by the choice of $\de$.   Finally, let $t\in \mc B$ and let us compute $\sum_{i\in t}(\mu(i))^{1/2}$:  First
of all we have that  $\sum_{i\in t}(\mu(i))^{1/2}=\sum_{i\in u}(\mu(i))^{1/2}$, where $u=t\cap P$. Now, since
$u\in \overline{\mc B\rest P}\con \overline{\mc D\otimes (n+1)}$, we can find $t_0<\dots <t_n$ in $\mc D$
such that $u\ip t_0\cup\dots\cup t_n$, and hence
\begin{equation}
\langle\mu^{1/2}, t\rangle=\sum_{j=0}^n\langle \mu^{1/2}, t_j\rangle \le
\frac{(n+1)(1+\vep)}{(nk)^\frac12},
\end{equation}
as promised. \fprue \cor\label{ewurehgree} Suppose that $\mc B$ is an $\al$-uniform barrier on $M$, $\al\ge
1$.  Then for every $\vep>0$ there is some $k=k(\al,\vep)$ such that  for every $\be$-uniform barrier on $M$
with $\be> \al k$ there $N\con M$ and some $U$-mean-assignment $\mu:\mc C\rest N\to c_{00}^+$ such that,
\begin{equation}
\sup\conj{ \langle \mu(s)^{1/2}, t\rangle}{t\in \mc B}\le \vep
\end{equation}
for every $s\in \mc B\rest N$. \qed \fcor

\lema\label{ertrffs} Fix an indecomposable countable $\al$ and a
sequence $(\vep_n)$ of positive reals. Then:

\noindent (a) there is a collection $(\mc B_n)$ of $\al_n$-uniform barriers on $\N/n$ and a corresponding
sequence of $U$-mean-assignments $\mu_n:\mc B_n \to c_{00}^+$ with the following properties:

\noindent\, (a.1) $\al_n>0$, $\sup_n\al_n=\al$,


\noindent \, (a.2) for every $m<n$ and every $s\in \mc B_n $
\begin{equation}
\sup\conj{\langle \mu_n(s)^{\frac12},t\rangle}{t\in \mc B_m}<\vep_n.
\end{equation}

\noindent (b)  Suppose that in addition $\al=\om^\ga$ with $\ga$ limit. Let $\al_n\uparrow \al$ be any
sequence such that $\al_n \om \le \al_{n+1}$ ($n\in \N$). Then there is a double sequence $(\mc B_i^n)$ such
that for every integers $n$ and $i$

\noindent \, (b.1) $\mc B_i^n$ is an $\al_i^{(n)}$-uniform barrier on $\N/(n+i)$, with $\al_i^{(n)}>0$ and
$\al_i^{(n)}\uparrow_i \al_n $.

\noindent \, (b.2) There are $U$-mean-assignments $\mu_{n,i}:\mc B_i^n\to c_{00}$ such that for every $s\in
\mc B_i^n$, and every $(m,j)<_\mathrm{lex} (n,i)$
\begin{equation}
\sup\conj{\langle \mu_{n,i}(s)^\frac12, t\rangle}{t\in  \mc B_j^m}<\vep_{n+i},
\end{equation}
where we recall that $<_\mathrm{lex}$ denotes the lexicographical order on $\N^{2}$ defined by
$(m,i)<_\mathrm{lex} (n,j)$ iff $m<n$, or $m=n$ and $i<j$.

\flema
\prue
(a): Choose  $\al_n \uparrow_n \al$ such that for every $n\in \N$, $\al_{n+1}> \al_{n}k(\al_{n},\vep_n)$,
that
 is is possible since $\al$ is indecomposable. Let $\mc C_n$ be an  $\al_n$-uniform family on $\N$ ($n\in
\N$). By Corollary \ref{ewurehgree} we can find a fusion sequence $(M_n)$ such that

\noindent (c) $\mc C_m\rest M_m \con \overline{\mc C_n}  $ if
$m\le n$, and

\noindent (d) for every $n\in \N$ there is a $U$-mean-assignment $\nu_n:\mc C_n\rest M_n\to c_{00}^+$ such
that
\begin{equation}
\sup\conj{\langle \nu_n(s)^{\frac12},t\rangle}{t\in \bigcup_{l<n}\mc C_{l}}<\vep_n
\end{equation}
for every $s\in \mc C_n\rest M_n$. Let $M=\{m_n\}$ be the fusion set of $(M_n)$, and $\Theta:M\to \N$ be the
corresponding order preserving onto mapping. It is not difficult to see that $\mc C_n=(\Theta"\mc B_n)\rest
(\N/n)$, and $\mu_n: \mc C_n\to c_{00}$ defined naturally out of $\nu_n$  $\Theta$ fulfils all the
requirements.

(b): Suppose that $\al=\om^\ga$ with $\ga$ limit. Let $\al_n\uparrow \al$ be any sequence such that
$\al_n \om \le \al_{n+1}$ ($n\in \N$).

\clam There is a fusion sequence $(M_n)$, $M_n=\{m_i^{(n)}\}$,  a
double sequence $(\mc B_i^n)$ of $\al_i^{(n)}$-uniform barriers on $M_n/m_i^{(n)}$ and $U$-mean-assignments
$\mu_{n,i}:\mc B_i^n\to c_{00}^+ $ such that

\noindent (e) $\al_i^{(n)}\uparrow_i\al_n$ ($n\in \N$), and

 \noindent (f) for every $(m,j)<_\mathrm{lex}(n,i)$, every $s\in \mc B_i^n$ and every
$t\in \mc B_j^m$, $\langle (\mu_{n,i}(s))^{1/2},t\rangle<\vep_{n+i}$.
\fclam
\prucl
First, use Corollary \ref{ewurehgree}  applied to $\al_0$ to produce an infinite set $M_0=\{m_i^{(0)}\}$ and
a sequence $(\mc B_i^0)$ of $\al_i^{(0)}$-uniform barriers on $M_0/\{m_i^{(0)}\}$ with $\al_i^{(0)}\uparrow
\al_0$ and $U$-mean-assignments $\mu_{0,i}:\mc B_i^0\to c_{00}$ such that for every $i$ and every $s\in \mc
B_i^{0}$, $\langle \mu_{0,i}(s)^{1/2},t\rangle \le \vep_i$ for every $t\in \mc B_j^0$ with $j<i$. In general,
suppose we have found  for every $k\le n$ $M_k=\{m_i^{(k)}\}\con M_{k-1}$, $(\mc B_i^k)$
$\al_i^{(k)}$-uniform barriers on $M_k/m_i^{(k)}$ and $U$-mean-assignments $\mu_{k,i}:\mc B_i^k\to c_{00}$
such that for every $(k,j)<_\mathrm{lex} (m,i)$  every $s\in \mc B_i^{m}$ and every  $t\in \mc B_j^k$
$\langle \mu_{m,i}(s)^{1/2},t\rangle \le \vep_{m+i}$. For each $k\le n$ define the following families
\begin{equation}
\mc B_k=\conj{s\con M_k}{{_*}s\in \mc B_{\min s}^k}.
\end{equation}
This is clearly an $\al_k$-uniform family on $M_k$. Since $\al_n\om \le \al_{n+1}$,  we can use again
Corollary \ref{ewurehgree} and find an infinite subset $M_{n+1}=\{m_i^{(n+1)}\}\con M_n$ and a sequence $(\mc
B_{i}^{n+1})$ of $\al_i^{(n+1)}$-uniform barriers on $M_{n+1}/m_{i}^{(n+1)}$ and $U$-mean-assignments
$\mu_{n+1,i}: \mc B_i^{n+1}\to c_{00}$ such that for every $s \in \mc B_i^{n+1}$,
\begin{equation}
\sup \conj{\langle (\mu_{n+1,i}(s))^{\frac12}, t\rangle}{t\in \bigcup_{k\le n} \mc B_{m} \cup
\bigcup_{j<i} \mc B_{j}^{(n+1)}}< \vep_{n+i+1},
\end{equation}
so, in particular for every $k\le n$ and every $t\in \mc B_{j}^{k}$, $\langle
(\mu_{n+1,i}(s))^{\frac12}, t\rangle<\vep_{n+i+1}$.  \fprucl
 Let $M$ be the fusion set of $(M_n)$, i.e. $M=\{m_0^{(n)}\}$. Observe that $m_0^{(n+i)}\ge m_{i}^{(n)}$ for
 every $n$ and $i$, so $M/m_0^{(n)}\con M_n/m_{i}^{(n)}$. Set $\mc C_i^{n}=\mc B_i^{n}\rest (M/m_0^{(n+i)})$.
This is an $\al_i^{(n)}$-uniform barrier on $M/m_0^{(n+i)}$. Consider $\nu_{n,i}=\mu_{n,i}\rest \mc
C_i^{n}:\mc C_i^{n}\to c_{00}$ has the property that for every $(m,j)<_\mathrm{lex} (n,i)$, every
every $s\in \mc C_i^n$ and every $t\in \mc C_j^m$, $\langle
(\nu_{n,i}(s))^{1/2},t\rangle<\vep_{n+i}$. Now use $\Theta:M \to \N$, $\Theta(m_0^{(n)})=0$, to
define the desired mean-assignments and  families.
\fprue

\nota
Observe that if $\mc B$ is $\al$-uniform on $M$ with $\al>0$, then $M^{[1]}\con \overline{\mc B}$. It readily
follows that the mean-assignments $\mu_n$ and $\mu_{n,i}$ obtained in Lemma \ref{ertrffs} have the property
that $\nrm{\mu_n(s)^{1/2}}_\infty \le \vep_n$ and $\nrm{\mu_{n,i}(s)^{1/2}}_\infty\le \vep_{n+i}$ for every
$s$ in the corresponding domains.
\fnota

\prop \label{dfjklweobf}

\noindent (a) Suppose that $\mc C$ and $\mc B_i$ are $\be$ and $\al_i$-uniform families on $M$
($i\in \N$) with $\al_i\uparrow \al$, $\al_i,\be \ge 1$. Let $\sig:\fin^{[<\infty]}\to \N$ be 1-1.
Then for every   $n\in \N$ the family
\begin{align*}
\mc D=\{& s_0\cup \dots\cup s_{n}\, :\, (s_i) \text{ is block, } s_0\in \mc C\text{ and } s_{i}\in
\mc B_{\sig((s_0,\dots,s_{i-1}))}\text{ for every $1\le i\le n-1$}\}
\end{align*}
is $\ga$-uniform on $M$, where $\ga=\al n+ \be^-$ if $1\le \be <\om$ and $n>0$, and $\ga=\al n +
\be$ if $\be\ge \om$ or $n=0$.

\noindent (b) Suppose that  $\mc B_i$ is $\al_i$-uniform  on $M$  ($i\in \N$) with $\al_i\uparrow
\al$. Let $\sig:\fin^{[<\infty]}\to \N$ be 1-1. Then  the family
\begin{align*}
\mc C=\{ \{n\}\cup s_0\cup \dots\cup s_{n-1}\, :& (\{n\},s_0,\dots,s_{n-1}) \text{ is block,  and }\\
& s_{i}\in \mc B_{\sig((\{n\},s_0,\dots,s_{i-1}))}\text{ for every $0\le i\le n-1$}\}
\end{align*}
is $\al\om$-uniform on $M$.

\fprop
\prue
(a): The proof is by induction on $n$. If $n=0$, the result is clear. So suppose that $n>0$. Now
the proof is by induction on $\be$.  Suppose first that $\be=1$. Then $\mc C=M^{[1]}$, and so,  for
every $m\in M$
\begin{align*}
\mc D_{\{m\}}=\{  s_1\cup \dots\cup s_{n}\, : & (s_1, s_2,\dots,s_n) \text{ is block, } s_1\in \mc
B_{\sig((\{m\}))}\text{ and }\\
&  s_{i}\in \mc B_{\sig((\{m\},s_1,s_2,\dots,s_{i-1}))}\text{ for every $2\le i\le n-1$}\},
\end{align*}
so, by inductive hypothesis, $\mc D_{\{m\}}$ is $\al (n-1)+\ga_m$-uniform on $M/m$, depending
whether $\al_m$ is finite or infinite, but in any case with $\ga_m\uparrow \al$. Hence $\mc D$ is
$\al n$-uniform on $M$. The general case for $1\le \be<\om$ is shown in the same way.

Suppose now that $\be\ge \om$. Then for every $m\in M$
\begin{align*}
\mc D_{\{m\}}=\{ & t\cup s_1\cup \dots\cup s_{n}\, :\, (t, s_1,\dots,s_n) \text{ is block, } t\in
\mc
C_{\{m\}}\text{ and }\\
&  s_{i}\in \mc B_{\sig((\{m\}\cup t,s_1,\dots,s_{i-1}))}\text{ for every $1\le i\le n-1$}\},
\end{align*}
By inductive hypothesis, $\mc D_{\{m\}}$ is $\al n+ \ga_m$-uniform on $M/m$, with $\ga_m \uparrow
\be$ , so $\mc D$ is $\al n+\be$-uniform on $M$, as desired.

(b) follows easily from (a).
\fprue

The following is a  generalization of Maurey-Rosenthal example for arbitrary countable indecomposable ordinal
$\al$.
\teor\label{lertjiorjgf}
For every countable indecomposable ordinal   $\al$ there is a normalized weakly-null sequence which is
$\be$-unconditionally saturated for every $\be<\al$ but without unconditional subsequences.
\fteor
\prue
Our example is a slightly modification of a $U$-sequence introduced in \cite{lop-tod}. So, we are going to
define a $\al$-uniform barrier $\mc B$ on $\N$, a $U$-mean-assignment $\vphi:\mc B\to c_{00}$  and some $\mc
G\con \fin\times \fin$ and then define the norm on $c_{00}$ by
\begin{equation}\label{ljeijgidjf}
\nrm{\xi}= \max\{\nrm{\xi}_{\infty},\sup\conj{|\langle \vphi(s)\rest t,\xi \rangle|}{(s,t)\in \mc G } \}
\end{equation}
where $\mc G\con \fin\times \fin$ is such that its first projection is $\mc B$. Notice that some sort of
restrictions have to be needed in the formula \eqref{ljeijgidjf}, since it is not difficult to see that
that for a compact and hereditary family $\mc F$, a normalized weakly-null sequence $(x_i)_i$ is $\mc
F$-unconditional iff it is equivalent to the evaluation mapping sequence $(p_i)_i$ of a weakly-compact subset
$K\con c_0$ that is $\mc F$-closed, i.e. closed under restriction on elements of $\mc F$.

Fix $\vep>0$, and let $\vep_n=\vep/2^{n+3}$.  Suppose that $\al=\om^\ga$. There are two cases to consider.
Suppose first that $\ga=\be+1$. We apply Lemma \ref{ertrffs} (a) to the indecomposable ordinal $\om^\be$ and
$(\vep_n)$ to produce the corresponding sequences of barriers $(\mc C_n)$ and  $U$-mean-assignments
$\nu_n:\mc C_{n}\to c_{00}$ ($n\in \N$) satisfying the conclusions (a.1) and (a.2) of the Lemma. If $\ga$ is
limit, then we use the part (b) of that lemma to produce a double sequence $(\mc B_i^n)$ and
$U$-mean-assignments $\nu_{n,i}:\mc C_{i}^n\to  c_{00}$ satisfying (b.1) and (b.2). In order to unify the two
cases we set for $n,i$,
\begin{equation*}
\mc B_{i}^n=\left\{\begin{array}{ll} \mc C_{i} & \text{if $\ga$ is successor ordinal} \\
\mc C_{i}^n & \text{if $\ga$ is  limit ordinal}
\end{array}\right.
\end{equation*}
and
\begin{equation*}
\mu_{n,i}=\left\{\begin{array}{ll} \nu_{i} & \text{if $\ga$ is successor ordinal} \\
\nu_{n,i} & \text{if $\ga$ is  limit ordinal}.
\end{array}\right.
\end{equation*}
Let $\sig:\fin^{[<\infty]}\to \N$ be  1-1 mapping such that $\sig((s_0,\dots,s_{n}))>\max s_n$ for
every block sequence $(s_0,\dots,s_n)$ of finite sets. For each $n$ define
\begin{align*}
\mc C_n=\{& s_0\cup\dots\cup s_{n-1}\, :\, (s_i)\text{ is block and  }\text{$s_i\in \mc
B_{\sig((\{n\},s_0,\dots,s_{i-1}))}^n$ for every $0\le i\le n$ }\},
\end{align*}
So, by Proposition \ref{dfjklweobf}, if $\al=\om^{\be+1}$, then $\mc C_n$ is a $\om^\be (n-1)+\zeta$-uniform
family  on $\N$, where $\zeta$ is such that $\mc B_{\sig((\{n\}))}^n$ is $\zeta$-uniform; while if
$\al=\om^\ga$ with $\ga$ limit, then it is $\al_n (n-1)+\zeta$ where $\zeta$ is such that $\mc
B_{\sig((\{n\}))}^n$ is $\zeta$-uniform.  Now  let
\begin{equation}
\mc C=\conj{s\in \fin}{{_*}s\in \mc C_{\min s}}.
\end{equation}
It turns out that  $\mc C$ is an $\al$-uniform family on $\N$ (so it is a front), not necessarily a barrier.
Observe that every $s\in \mc C$ has a unique decomposition $s=\{n\}\cup s(0)\cup \dots\cup s(n-1) $  with
$n=\min s$ and  $s(i)\in \mc B_{\sig(s[i])}$, and where $s[i]=(\{n\},s_0,\dots,s_{i-1})$  ($0\le i\le n-1$).
For every $s\in \mc C$ and every $i\le s$, set
$$\xi(s,i)=(\mu_{\min s, \sig(s[i])}(s(i)))^{1/2}.$$
Define now $\Phi:\mc C\to c_{00}$ for  every $s\in \mc C$   by
\begin{equation}
\Phi(s)=e_{\min s}+ \sum_{i=0}^{n-1}\xi(s,i),
\end{equation}
It is not difficult to see that $\Phi:\mc C\to c_{00}$ is a $U$-mapping. Now define on $c_{00}$ the norm
\begin{align}\label{pjreigjhjgf}
\nrm{\xi}=&\sup\conj{|\langle  \Phi(s)\rest (s\setminus t),\xi \rangle|}{s\in \mc C,   t\con s(i), \text{ for some $i<\min s$}}=\nonumber\\
=&\sup\conj{|\langle  \Phi(s)\rest (u\setminus t),\xi \rangle|}{u\ip s\in \mc C,   t\con s(i), \text{ for some $i<\min s$}},
\end{align}
the last equality because $\Phi$ is Lipschitz and supported by a front.
 Let $\eqs$ the completion of $c_{00}$ under this norm. Then the Hamel basis $(e_n)_n$ of $c_{00}$ is a
normalized basis of $\eqs$, moreover monotone (since $\Phi$ is Lipschitz with domain a front) and
weakly-null: To prove this, it is enough to  see that the set
$$L=\conj{\Phi(s)\rest(u\setminus t)}{s\in \mc B,\, u\ip s,\text{ and $t\con s(i)$ for some $i<\min s$}}$$
is weakly-compact. So, let $(\Phi(s_n)\rest (u_n\setminus t_n))_n$ a typical sequence in $L$. Since $\mc C$
is a front, we can find an infinite set $M$ and  $u\in \fin $ such that $(u_n)_{n\in M}$ converges to $u$ and
such that $(s_n)$ is a $\De$-system with root $u\ip r$. Since $\Phi$ is Lipschitz de, we obtain
 that $(\Phi(s_n)\rest t_n)_{n\in M}$ converges to $\Phi(s_m)\rest t$ for (any) $m\in M$. If $u=\buit$, then $(\Phi (s_n)\rest (t_n\cup
 u))_{n\in M}$ converges to $0$. Otherwise,  let $N\con M$  and  $j<\min u$ be such that    $t_n \con
 s_n(j)$ for every $n\in N$. Now $(t_n)_{n\in N}$ is a sequence in the closure of $\mc B_{\sig(s[i])}^{\min s}$, hence, we can find  $P\con N$
such that $(t_n)_{n\in P}$ is convergent with limit $t$.  It follows that $(\Phi(s_n)\rest (u_n\setminus
t_n))_{n\in P}$ has limit $\Phi(s_n)\rest (u\setminus t)\in L$, where $n$ is (any) integer in $P$.

The next is a crucial computation.
\clam  For every $s,t\in \mc C $ and every $i\le \min s$ and $j\le \min s$,  we have that
\begin{equation*}
0\le \langle \xi(t,j),\xi(s,i)\rangle\le \left\{\begin{array}{ll} \vep_{\max\{\min s,\min t\}} & \text{if
$t[j]\neq s[i]$}\\
1 & \text{if $t[j]=s[i]$}.
\end{array}\right.
\end{equation*}
\fclam
\prucl Set $n=\min s$, $m=\min t$, and assume that $t[j]\neq s[i]$.  Suppose first that $\al=\om^{\be+1}$.
Then, by definition of the mean assignments, $\langle \xi(t,j),\xi(s,i)\rangle\le
\vep_{\max\{\sig(t[j]),\sig(s[i])\}}$, but $\sig(u_0,\dots,u_k)\ge \max u_k$ for every block sequence
$(u_i)$, which derives into the desired inequality. Assume now that $\al=\om^{\ga}$, $\ga$ limit ordinal. If
$\min s=\min t$, then $\langle \xi(t,j),\xi(s,i)\rangle\le \vep_{\min s+\max\{\sig(t[j]),\sig(s[i])\}}\le
\vep_{\min s}$.  While if $\min t\neq \min s$, say $\min t<\min s$, then $\langle \xi(t,j),\xi(s,i)\rangle\le
\vep_{\min s+\sig(s[i])}\le \vep_{\min s}$.

If $\sig(s[i])=\sig(t[j])=l$, then $\min s=\min t=n$, and
\begin{equation}
\langle \xi(s,i),\xi(t,j)\rangle\le \nrm{(\mu_{n,l}(s(i)))^{1/2}}_{\ell_2}
\nrm{(\mu_{n,l}(t(j)))^{1/2}}_{\ell_2}\le 1,
\end{equation}
since both are means.
\fprucl
\clam
The summing basis  $(S_n)$ of $c$ is finitely block represented in any subsequence of $(e_n)_n$.
\fclam
\prucl
Fix an infinite set $M$ of integers, and $l\in \N$. Let   $v\in \mc B\rest M/l$, $v=\{n\}\cup v(0)\cup \dots
\cup v(n-1)$ its canonical decomposition, and set
\begin{equation}
x_i=\sum_{j\in v(i)} \xi(v,i)(j)e_j.
\end{equation}
Observe that $\langle \Phi(v),x(v,i) \rangle=  \langle  \xi(v,i),\xi(v,i) \rangle= 1$, so from the previous
claim we obtain that $\nrm{x_i}=1$. Now consider scalars $(a_i)_{i\le n-1}$ with $\nrm{\sum_{i\le n-1}a_i
S_i}_{\infty}=1$. Observe that this implies that $\max_{i\le n-1}|a_i|\le 2$.   We are going to show that
\begin{equation}
1 \le \nrm{\sum_{0\le i\le n-1} a_i x_i} \le
3+\vep.
\end{equation}
To get the left hand inequality, suppose that $1=\nrm{\sum_{0\le i\le n-1} a_i S_i}_\infty=|\sum_{i\le m
}a_i|$, where $m\le n-1$. Let $t=\{n\}\cup s(0)\cup \dots \cup s(m)$. By \eqref{pjreigjhjgf} it follows that
\begin{equation}
\nrm{ \sum_{i\le n-1}a_i x_i}\ge \langle \Phi(v)\rest t, \sum_{i\le n-1}a_i x_i \rangle=|\sum_{i\le m}a_i|=1.
\end{equation}
Next,  fix $s\in \mc C$  and $t\con s(k)$ for some $k<\min s$.  Suppose first that $\min v=\min s$. Let
$i_0=\max\conj{i\le n-1}{v(i)=s(i)}$. If $k>i_0$  then by the previous claim  we obtain
\begin{align} |\langle \Phi(s)\rest (s\setminus t), \sum_{i\le n-1}a_i
x_i \rangle)|\le & |\sum_{i\le i_0 }a_i|+ \sum_{i_0<i,j\le n-1} 2|\langle \xi(s,i),
\xi(t,j)|\rangle \le \nonumber\\ \le & \nrm{\sum_{i\le n-1}a_i S_i}_{\infty} +  2 n^2\vep_n\le
  (1+\vep)\nrm{\sum_{i\le n-1}a_i S_i}_{\infty}.
\end{align}
Suppose that $k\le i_0$. Then
\begin{align} |\langle \Phi(s)\rest (s\setminus t), \sum_{i\le n-1}a_i
x_i \rangle)|\le & |\sum_{i\le i_0, i\neq k }a_i+ a_k \langle \xi(v,k),\xi(v,k)\rest (s(k)\setminus t) \rangle|+ \nonumber \\
& +  \sum_{i_0<i,j\le n-1} 2|\langle
\xi(s,i),
\xi(t,j)|\rangle \le \nonumber\\ \le & 3\nrm{\sum_{i\le n-1}a_i S_i}_{\infty} +  2 n^2\vep_n\le
  (3+\vep)\nrm{\sum_{i\le n-1}a_i S_i}_{\infty}.
\end{align}
Suppose now that  $n=\min v\neq \min s$, say $\min s<\min v$.  Let $i_0<n$, if possible, be such that $\min
s\in v(i_0)$. Then,
\begin{align}
|\langle \Phi(s)\rest (s \setminus t), \sum_{i\le n-1}a_i
x_i \rangle)|\le    &   |a_{i_0}|\nrm{\xi(v,i_0)}_\infty +2 \sum_{i_0\le i<n}\sum_{0\le j<\min
t}\langle\xi(t,j),\xi(s,i)\rangle \le \nonumber \\ \le &  2 \vep_n+2 n^2\vep_n\le  \vep.
\end{align} \fprucl
Finally, it rests to show that the sequence  $(e_n)$ is $\be$-unconditionally saturated for every $\be<\al$.
We consider the two obvious cases:

\noindent \textsc{Case 1.} $\al=\om^{\be+1}$.  Let
$$\mc
D=\conj{s\con \N}{{_*}s\in \mc B_{\min s}^{0}}.$$
This is  an
$\om^\be$-uniform family on $\N$ since each family $\mc B_m^0$  is $\al_m$-uniform and   $\sup_m \al_m =\om^\be$. Therefore, the  next claim  gives that $(e_n)$
is $\be$-unconditionally saturated for every $\be<\al$.

\clam    $(e_n)_n$   is  $\mc D$-unconditional with
constant at most $2+\vep$.
\fclam
\prucl  Fix $t\in \mc D$, and
let $(a_i)_{i\in \N}$ be scalars such that $\nrm{\sum_{i\in \N} a_i e_i} =1$. Fix also $s\in \mc C$. Suppose
first that $\min s \in t$. Then since $\sig(s[i])>\min s\ge \min t$ and ${_*}t\in \mc B_{\min t}^{0}$ we
obtain that
\begin{equation}
|\langle \Phi(s),  \sum_{i\in t}a_i e_i\rangle|  \le |a_{\min s}|+\vep \le (1+\vep)\nrm{\sum_{i}a_i e_i}.
\end{equation}
Now suppose that $\min s \notin t$, but    $s\cap t\neq \buit$ (otherwise $\langle \Phi(s),  \sum_{i\in t}a_i
e_i\rangle=0$). Let
$$i_0=\min \conj{i< \min s}{ s(i)\cap t\neq \buit}.$$ Then for every $i_0<i<\min s$ we have that
$\sig(s[i])>\max s_{i_0}\ge \min t$, so
\begin{equation}
|\sum_{j\in t}a_j \xi(s,i)(j)|<\vep_{\sig(s[i])},
\end{equation}
hence
\begin{align}
|\langle \Phi(s)\rest u,  \sum_{i\in t}a_i e_i\rangle|\le & |\sum_{j\in t\cap s(i_0)} a_j \xi(s,i_0)(j)|+\sum_{i_0<i<\min s}|\sum_{j\in
t}a_j \xi(s,i)(j)| = \nonumber \\
=& |\langle \Phi(s)\rest (\{n\}\cup s(0)\cup \dots \cup (s(i_0)\cap t)),\sum_{i\ge \min t}a_ie_i \rangle|   + \nonumber \\
+ & \sum_{i_0<i<\min s}|\sum_{j\in
t}a_j \xi(s,i)(j)|
\le   \nrm{\sum_{i\ge \min t}a_i e_i} + \vep \nrm{\sum_{i\in \N}a_ie_i} \le \nonumber \\
\le & (2 +\vep) \nrm{\sum_i a_i e_i},
\end{align}
the last inequality because $(e_i)$ is monotone.
\fprucl
\noindent \textsc{Case 2.} $\al=\om^\ga$, $\ga$ a countable limit
ordinal. The desired result follows from the following fact.
\clam
For every $n\in \N,$ the sequence  $(e_i)$ is $\mc B_{0}^n$-unconditional with constant at most $2n+1$.
\fclam
\prucl Fix $n\in \N$ and  $t\in \mc B_0^n$. Let
$(a_i)_{i\in \N}$ be scalars such that $\nrm{\sum_{i\in \N} a_i e_i}=1$. Fix  $s\in \mc C$. Suppose first
that $n\le \min s$. Then in a similar manner that in \textsc{Case 1} one can show that
\begin{equation}
|\langle \Phi(s),  \sum_{i\in t}a_i e_i\rangle|\le |a_{\min s}|+\vep \le (1+\vep)\nrm{\sum_{i}a_i e_i}.
\end{equation}
Suppose that $m=\min s < n$, then
\begin{align}
|\langle \Phi(s),  \sum_{i\in t}a_i e_i\rangle|\le &  |a_{\min s}|+\sum_{i=0}^{m-1} |\sum_{j\in s(i)\cap t } a_j \xi(s,i)(j)|= \nonumber \\
 = &  |a_{\min s}|+\sum_{i=0}^{m-1} |\langle \Phi(s)\rest u_i, \sum_{j\ge \min (s(i)\cap t)}a_j  \rangle|   \le \nonumber \\
 \le  & (2m +1)\nrm{\sum_i a_i e_i}.
\end{align}
where $u_i= s(0)\cup \dots \cup(s(i)\cap t)$.
\fprucl
\fprue

\cor
For every indecomposable ordinal   $\al$ there is a weakly-compact $K\con c_{00}$ such that

\noindent (a) $K\con B_{c_0}$ is point-finite (i.e. $\conj{\xi(n)}{\xi \in K}$ is finite for every integer
$n$) supported by a $\al$-uniform barrier on $\N$,

\noindent (b) the evaluation mapping sequence $(p_i)_i$ of $C(K)$ is a normalized weakly-null monotone basic
sequence, and

\noindent (c) The summing basis of $c$ is $4$-finitely representable in every subsequence of $(p_i)_i$;
hence no subsequence of $(p_n)$  is unconditional, but

\noindent (d) $(p_i)_i$ is $\be$-unconditionally saturated for every $\be<\al$.
\fcor
\prue
Let $\mc C$  be  the $\al$-uniform family on $\N$ and let $\Phi:\mc C\to c_{00}$ be the $U$-mapping given in
proof of Theorem \ref{lertjiorjgf}. Let $M\con \N$ be such that $\mc C\rest N$ is a $\al$-uniform barrier on
$\N$. Let $\theta $ be the order-preserving mapping from $M$ onto $\N$. Let $\mc B=\theta"\mc
C=\conj{\theta"(s)}{s\in \mc C}$ and let $\vphi: \mc B\to c_{00}$ be naturally defined by
$\vphi(s)=\Phi(\theta^{-1}(s))$.  $\mc B$ is a uniform barrier on $\N$ and $\vphi$ is a $U$-mapping. Observe
that every $s\in \mc B$ has a unique decomposition, given by the one of $\theta^{-1}s$.        Let
$$K=\conj{\vphi(s)\rest (u\setminus t)}{u\ip s\in \mc B,\, t\con s(i)\text{ for some $i$}}.$$
This is a weakly-compact subset of $c_0$, and the corresponding evaluation mapping sequence $(p_i)_i$  is
1-equivalent to the subsequence $(e_n)_{n\in M}$ of the weakly-null sequence $(e_i)_i$ given in the proof of Theorem \ref{lertjiorjgf}. So $K$
fulfills all the requirements.
\fprue


\begin{thebibliography}{this}
\bibitem{aa}
D. Alspach and S.A. Argyros, {\em Complexity of weakly null sequences}, Dissertationes
Mathematicae, {\bf 321}, (1992), 1--44.



\bibitem{arg-mer-tsa} S.A. Argyros, S. Mercourakis and A. Tsarpalias,  \emph{Convex unconditionality and summability
of weakly null sequences}, Israel J. Math. \textbf{107} (1998), 157--193




\bibitem{arg-god-ros} S. A. Argyros; G. Godefroy; H. P. Rosenthal,
\emph{Descriptive set theory and Banach spaces}.   Handbook of the geometry of Banach
spaces, Vol. 2, 1007--1069, North-Holland, Amsterdam, 2003.

\bibitem{tod1}S. A. Argyros and S. Todorcevic, \emph{Ramsey methods in analysis}.  Birkh\"{a}user Verlag, Basel 2005.

\bibitem{arg-ka1}S. A. Argyros and V. Kanellopoulos, \emph{Determining $c_0$ in $C(K)$ spaces},  Fund. Math.  187  (2005),  no. 1, 61--93.


\bibitem{arv} A. Arvanitakis, \emph{Weakly null sequences with an unconditional subsequence.} Preprint 2004.


\bibitem{bess-pel} C. Bessaga and A. Pe\l czy\'nski,  \emph{On bases and unconditional convergence of series in Banach spaces}. Studia Math. \textbf{17} 1958 151--164.

\bibitem{bess-pel-1} C. Bessaga and A. Pe\l czy\'nski,  \emph{Spaces of continuous functions.
 IV. On isomorphical classification of spaces of continuous functions.} Studia Math. \textbf{19} 1960 53--62.

\bibitem{cas-shu} P.G. Casazza and T.J. Shura, \emph{Tsirelson's space}. Lecture Notes in Math., Vol.1363,
Springer-Verlag, Berlin 1989.


\bibitem{elton} J. Elton, \emph{Thesis}, Yale University (1978).


\bibitem{gas} I, Gasparis, \emph{ A dichotomy theorem for subsets of the power set of the natural
numbers.} Proc. Amer. Math. Soc.  \textbf{129}  (2001),  no. 3, 759--764.
\bibitem{ga-od-wa} I. Gasparis, E. Odell and B. Wahl, \emph{Weakly null sequences in the Banach space $C(K)$}. Preprint 2004.




\bibitem{li-tza}J. Lindenstrauss and L. Tzafriri,
\emph{ Classical Banach spaces. I. Sequence spaces}. Ergebnisse der Mathematik und ihrer
Grenzgebiete, Vol. 92.

\bibitem{lop-tod} J. Lopez Abad and S. Todorcevic,  \emph{Partial unconditionality of weakly null sequences}.  RACSAM Rev. R. Acad. Cienc. Exactas F\'{\i}s. Nat. Ser. A Mat.  100  (2006),  no. 1-2.
\bibitem{nashwill} C. St. J. A. Nash-Williams, \emph{On well-quasi-ordering transfinite sequences,}
 {Proc. Cambridge Philos. Soc.} {\bf 61} (1965), 33--39.
\bibitem{mau-ros}B. Maurey and  H. P. Rosenthal,
\emph{Normalized weakly null sequence with no unconditional subsequence}. Studia Math. 61
(1977), no. 1, 77--98.

\bibitem{mau1} B. Maurey, \emph{Une suite faiblement convergente vers zero sans sous-suite
inconditionnelle}, S\'{e}minaire Maurey-Schwartz 1975-76, expos\'{e} IX (1976).

\bibitem{mau2} B. Maurey, \emph{Quelques resultats concernant l'inconditionnalite}, S\'{e}minaire Maurey-Schwartz 1975-76, expos\'{e} XVI (1976).



\bibitem{odell} E. Odell, \emph{Applications of Ramsey theorems to Banach space theory}, Notes in Banach spaces,
pp. 379--404, Univ. Texas Press, Austin, Tex., 1980.

\bibitem{odell1} E. Odell, On Schreier unconditional sequences, Contemp. Math. \textbf{144} (1993), 197-201.

\bibitem{pel-sem} A. Pe\l czy\'nski and Z. Semadeni,
\emph{Spaces of continuous functions III. Spaces $C(\Omega)$ for $\Omega$ without perfect subsets}, Studia
Math. \textbf{18} (1959), 211--222.


\bibitem{pud-rod}P. Pudlak and V. R\"{o}dl, \textit{Partition theorems
for systems of finite subsets of integers}, Discrete Math.,
\textbf{39}, (1982), 67-73.



\bibitem{rosenth} H. P. Rosenthal, \emph{The Banach spaces $C(K)$}.  Handbook of the geometry of Banach spaces, Vol. 2,  1547--1602, North-Holland, Amsterdam, 2003.




%
%
%
\end{thebibliography}
\end{document}